 \documentclass[11pt]{article}
 \usepackage{amsmath}
 \usepackage{pifont}
 \usepackage{amsfonts}
 \usepackage{mathrsfs}

 \usepackage{booktabs}
 \usepackage{threeparttable}
 \usepackage{mathrsfs,amsfonts,amsmath}
 \usepackage[dvips]{color}

\newcommand\email[1]{\href{mailto:#1}{ \nolinkurl{#1}}}

 \newtheorem{theorem}{Theorem}[section]
 \newtheorem{definition}[theorem]{Definition}
 \newtheorem{lemma}[theorem]{Lemma}
 \newtheorem{corollary}[theorem]{Corollary}
 \newtheorem{proposition}[theorem]{Proposition}
 \newtheorem{remark}[theorem]{Remark}
 \newtheorem{condition}[theorem]{Condition}
 \newtheorem{example}{Example}[section]

 \def\blemma{\begin{lemma}}\def\elemma{\end{lemma}}
 \def\bproposition{\begin{proposition}}\def\eproposition{\end{proposition}}
 \def\ttheorem{\begin{theorem}}\def\etheorem{\end{theorem}}
 \def\bcorollary{\begin{corollary}}\def\ecorollary{\end{corollary}}
 \def\bremark{\begin{remark}}\def\eremark{\end{remark}}
 \def\bcondition{\begin{condition}}\def\econdition{\end{condition}}

 \def\benumerate{\begin{enumerate}}\def\eenumerate{\end{enumerate}}
 \def\bitemize{\begin{itemize}}\def\eitemize{\end{itemize}}

 \def\beqlb{\begin{eqnarray}}\def\eeqlb{\end{eqnarray}}
 \def\beqnn{\begin{eqnarray*}}\def\eeqnn{\end{eqnarray*}}
 \def\ar{\!\!\!&}

 \def\mbb{\mathbb}

 \def\proof{\noindent{\it Proof.~~}}\def\qed{\hfill$\Box$\medskip}

\begin{document}

\noindent{Version: 2016/03/08}

\bigskip\bigskip

\centerline{\Large\bf Nonparametric Estimation for Jump-Diffusion CIR Model}

%

\bigskip

\centerline{Wei Xu}

\medskip

\centerline{School of Mathematical Sciences, Beijing Normal University}

\medskip

\centerline{Beijing 100875, People's Republic of China}

\medskip
\centerline{E-mails: \texttt{xuwei@mail.bnu.edu.cn}}
\bigskip

{\narrower

\noindent\textit{Abstract:} We study the nonparametric estimation for the
intensity of Poisson random measure in jump-diffusion CIR model based on the low frequency observations.  This is given in terms of the minimization of norms on a nonempty, closed and
convex subset of some special Hilbert space. We establish the measurability
of the estimator and derive its consistency and asymptotic risk bound.

\smallskip


\medskip


\par}

\section{Introduction and main results}

\setcounter{equation}{0}
\medskip

  The \textit{Cox-Ingersoll-Ross model} (CIR model) defined by the following stochastic differential equation (SDE) was firstly introduced by Cox et
  al. (1985) in the study of term structure of interest rate:
  \beqlb\label{NPCB1.1}
  dY(t) = b(\beta-Y(t))dt + \sqrt{2cY(t)} dB(t),
 \eeqlb
  where $\beta,b,c>0$ are given constants and $\{B(t): t\ge 0\}$ is a standard Brownian motion. Motivated by the study of jump risks which can not be ignored in the pricing of assets, Ahn and Thompson (1988) studied the following jump-diffusion process by adding a jump component into (\ref{NPCB1.1}):
  \beqlb\label{NPCB1.2}
  dX(t) = b(\beta-X(t))dt + \sqrt{2cX(t)} dB(t)+\int_0^{\infty}zN(dt, dz),
  \eeqlb
  where $N(dt,dz)$ is a Poisson random measure on $(0,\infty)^2$ with L\'evy intensity $n(dz)$ and $(1\wedge z)n(dz)$ is a finite measure on $(0,\infty)$; more details can be seen Duffie et al. (2000, 2003). In this paper we always assume $n(dz)$ is
  absolutely continuous with respect to Lebesgue measure, i.e. there exits a non-negative function $k(z)$
  satisfying $n(dz)=k(z)dz$.

  Actually, by Theorem~2.5 in Fu and Li (2010), the nonnegative solution $\{X(t):t\geq 0\}$ to (\ref{NPCB1.2}) uniquely exists and is a continuous-sate branching process with immigration (CBI processes) with transition semigroup $(P_t)_{t\geq 0}$ given by
  \beqlb\label{NPCB1.3}
\int_0^\infty e^{-\lambda y}P_t(x,dy)
 =
\exp\Big\{-xv_t(\lambda) - \int_0^t\psi(v_s(\lambda))ds\Big\},
 \eeqlb
 where  \beqnn
v_t(\lambda)=\frac{be^{-bt}\lambda}{b+c\lambda(1-e^{-bt})}
 \eeqnn
 and
 \beqlb\label{NPCB1.4}
\psi(z) = b\beta z + \int_0^\infty(1-e^{-zu})n(du).
 \eeqlb
 Otherwise, from Theorem~3.20 in Li (2011, p66) we have the semigroup $(P_t)_{t\geq 0}$ is ergodic, i.e. for any $x\geq 0$, $P_t(x,\cdot)$ converges to a probability measure $\eta$ on $[0,\infty)$ as $t\rightarrow \infty$ and the Laplace transform of $\eta$ is given by
 \beqlb\label{NPCB1.5}
 L_{\eta}(\lambda)=\exp\Big\{ -\int_0^{\infty}\psi(v_s(\lambda))ds\Big\}.
 \eeqlb
 For any finite set $\{t_1<
t_2< \cdots< t_n\}\subset \mbb{R}$ we can define the probability measure
$\mathbf{Q}_{t_1,t_2,\cdots,t_n}$ on $\mbb{R}_+^n$ by
 \beqlb\label{NPCB1.6}
 \mathbf{Q}_{t_1,t_2,\cdots,t_n}(dx_1,dx_2,\cdots,dx_n) = \eta(dx_1)P_{t_2-t_1}(x_1,dx_2)\cdots P_{t_n-t_{n-1}}(x_{n-1},dx_n).
 \eeqlb
  Since $\{\mathbf{Q}_{t_1,t_2,\cdots,t_n}: t_1< t_2< \cdots<t_n\in \mbb{R}\}$ is a consistent family, there is a stationary Markov process $\{Z(t): t\in \mbb{R}\}$ with finite-dimensional distributions given by (\ref{NPCB1.6}) and one-dimensional marginal distribution $\eta$. From Remark~2.6 in Li and Ma (2015) and Birkhoff's ergodic theorem, we have $\{Z(t): t\in \mbb{R}\}$ is ergodic. By a fairly simple (continuous time) coupling argument, with out loss of generality we always assume $X(t)$ defined by (\ref{NPCB1.2}) is a stationary and ergodic process.

  Before applying (\ref{NPCB1.2}) into practical problems, the key preparation is estimating $(\beta,b,c)$ and $n(dz)$. Since estimations for $(\beta,b,c)$ have been given by Huang et
  al. (2011), we just need to found some suitable estimations of $n(dz)$ with other parameters known.
  There are a lot of works about parameter estimations for the standard CIR-model and  a review had been given in Xu (2014) including
   the \textit{conditional least squares estimators} (CLSEs) and the \textit{maximum likelihood estimators} (MLEs) given by Overbeck and Ryd\'{e}n (1997) and Overbeck (1998). Here we only give a summary of some known works about nonparametric estimation for jump-diffusion processes.
  Unfortunately, limited works have been done in the nonparametric estimation in jump-diffusion CIR models compared with L\'{e}vy processes.  Watteel and Kulperger (2003) proposed and implemented an approach for
  estimating the jump distribution of the L\'{e}vy processes by fixed spectral cut-off
  procedure. The penalized projection method was applied in Figueroa-L\'{o}pez and Houdr\'{e} (2006) to estimate the L\'{e}vy density on a compact interval separated from the origin,
  based on a continuous time observation of the sample path throughout a time interval $[0,T]$.
  Moreover, Figueroa-L\'{o}pez (2009) used the projection method for discrete
  observations and provided minimum risks of estimation for smooth L\'{e}vy densities,
  as well as estimated on a compact interval separated from the origin. Comte and Genon-Catalot used a Fourier approach to construct an adaptive nonparametric estimators and
  provide bounds for the global $\mathbb{L}^2$-
  risk with high frequency data and low frequency data respectively; see Comte and Genon-Catalot (2009, 2010). Neumann and Reiss (2009) studied the nonparametric estimation for L\'{e}vy processes based on the empirical characteristic function. Jongbloed et al. (2005) considered a related low-frequency problem for the canonical function in Ornstein-Uhlenbeck processes driven by L\'{e}vy processes and a consistent estimator has been constructed.

  In this work, based on the low frequency observations at equidistant time
  points $\{k\Delta :k=0,1,\dots\}$ of a single realization, we establish some nonparametric estimators for the L\'evy density $n(dz)$ by minimazing the norms of the elements of a closed and convex subset in some special Hilbert space.
   For simplicity, we take $\Delta=1$ and denote the observation be $\{X_k:k=0,1,\cdots\}$
  but all the results presented below can be extended to the general case. We always assume all functions below are defined on $\mathbb{R}_+$.
  Let $\mu(dz)=(1\wedge z)dz$ and
   \beqlb\label{NPCB1.7}
   \mathcal{L}(\mu):=\{f(z):\mu(|f|):=\int_0^{\infty}|f(z)|\mu(dz)<\infty\}.
   \eeqlb
  We say
  a real-valued function  $f\in \mathscr{V}_{b}$, if for any $-\infty<a<b<\infty$, $$V_a^b(f):=\sup_{n\in\mathbb{N}}\sup_{p\in \mathcal {P}_n}\sum_{i=0}^{n-1}|f(x_{i+1}-f{x_i})|<\infty,$$
  where $\mathcal {P}_n=\big\{p=\{x_0,\dots,x_{n}\}: a=x_0< \cdots<x_n=b  \big\}$.
 We define the following convex subset of $\mathcal {L}(\mu)$
 \beqlb\label{NPCB1.8}
 K:=\{k(z)\in \mathcal {L}(\mu)\cap \mathscr{V}_{b} :k(z)\mbox{ is nonnegative and right-continuous} \}.
 \eeqlb
 Otherwise, a linear operator $\mathbf{T}$ is defined by
 \beqlb\label{NPCB1.9}
 \mathbf{T}f(\lambda):= \int_0^1\int_0^{\infty}(1-e^{-zv_s(\lambda)})f(z)dzds,\quad f(z)\in \mathcal{L}(\mu).
 \eeqlb
  Let
 \beqlb\label{NPCB1.10}
  L_{n}(\lambda)=\frac{1}{n}\sum_{k=1}^ne^{-\lambda X_k+X_{k-1}v(\lambda)+D},
 \eeqlb
 where $v(\lambda)=v_1(\lambda)$ and $ D=\frac{b\beta}{c}\log(1+c\lambda (1-e^{-b})/b)$.
 By the erogidicity of $X(t)$ and the continuous mapping theorem, for any $\lambda>0$ we have
 \beqlb\label{NPCB1.11}
 g_{n}(\lambda):=-\ln(L_{n}(\lambda))\xrightarrow[\rm d]{\rm a.s.} \mathbf{T}k(\lambda).
 \eeqlb
 Let $\{k_R\geq 0:R=1,2,\cdots\}$ be a increasing sequence of functions satisfying that $\mu(|\tilde{k}_R|)\leq R$. For any $R\geq 1$ define
 $$K_R:=\{f\in K : f\leq k_R\mbox{ and } {\rm V}_{i}^{i+1}(f)\vee{\rm V}_{(i+1)^{-1}}^{i^{-1}}(f)\leq R \mbox{ for any }i=1,2,\cdots\}.$$
 Here for any fixed $R\geq 1$, we establish the following estimator for $k(z)$:
 \beqlb\label{NPCB1.12}
 \hat{k}_{R,n}(z)=\arg\min_{f\in K_R} \int_0^{\infty}|g_{n}(\lambda)-\mathbf{T}f(\lambda)|^2w(\lambda)d\lambda,
 \eeqlb
 where $w(\lambda)$ is a bounded and non-negative weighted function with compact support, denote by $S_w$, and there exist $0\leq a< b<\infty$
 such that $[a,b]\subset S_w$. We give the main results in the following three theorems.
  \begin{theorem}\label{NPCBt1.1}
 For any $R\geq 1$, we have $\hat{k}_{R,n}(z)$ is well defined, i.e. $\hat{k}_{R,n}(z)$ exists uniquely and is measurable.
 \end{theorem}
 \begin{theorem}\label{NPCBt1.2}
 If $k(\cdot)\in K_R$ for some $R>0$, then the estimator $\hat{k}_{R,n}$ is strongly consistent. In details,
 \beqlb\label{NPCB1.13}
 \mu(|\hat{k}_{R,n}-k|)\overset{\rm a.s.}\longrightarrow 0.
 \eeqlb
 \end{theorem}
 \begin{theorem}\label{NPCBt1.3}
 Suppose $z^2n(dz)$ is a finite measure and $\exp\{-\int_0^{\infty}\psi(v_s(\lambda))ds\}<\infty$ for any $\lambda\in (-\frac{be^{-b}}{c(1-e^{-b})},\infty)$. If $k(\cdot)\in K_R$ for some $R>0$, then there exists a constant $C>0$ such that for $n$ large enough have
 \beqlb\label{NPCB1.14}
 \sqrt{n}\textbf{E}[\mu(|\hat{k}_{R,n}-k|)]< C.
 \eeqlb
 \end{theorem}
 \begin{remark}
 Conditions in Theorem~\ref{NPCBt1.3} can be weakened; i.e., if $\int_0^{\infty}\psi(v_s(v(2a)-2v(a)))ds<\infty$ for some $a>0$, we can choose a weighted function $w(\lambda)$
 with $S_w\subset [0,A]$.
 \end{remark}

 This paper is organized as follows. In Section~2, we will prove Theorem~\ref{NPCBt1.1}. The consistency and asymptotic risk bound of estimator (Theorem~\ref{NPCBt1.2} and \ref{NPCBt1.3}) will be proved in Section~3. 

 {\bf Notation:} In this paper, we denote $\mathbb{R}_+=[0,\infty)$ and $\mathbb{Q}$ be the set of all retional number. Moreover, $L_{Q}(\lambda)$ denotes the Laplace transform of the measure $Q$, $\overset{\rm
 a.s.}\rightarrow$ and $\overset{\rm d}\rightarrow$ mean converge almost surely and in distribution
 respectively. Similarly, $\overset{\rm
 a.s.}=$ and $\overset{\rm d}=$ mean equal almost surely and in distribution.

 \section{Existence, uniqueness and measurability}
 \setcounter{equation}{0}
 \medskip

 In this section, we will prove Theorem~\ref{NPCBt1.1} by identifying the estimator defined by (\ref{NPCB1.12}) exists uniquely and is measurable. Firstly, we recall a conclusion which can be found in many books about functional analysis.

 \begin{lemma}\label{NPCBt2.1}
If $S$ is a Banach space with norm $||\cdot||$, $M$ is a nonempty, closed, convex subset of $S$, then $M$ contains a unique element of smallest norm.
 \end{lemma}
 With this lemma we will give the most important theorem, which will guarantee the measurability of the estimators. Actually, least squares estimators (LSEs) and maximum likelihood
estimators (MLEs) are just special cases of this theorem.

 \begin{theorem}\label{NPCBt2.2}
 Suppose $(\Omega, \mathscr{F},\mathbb{P})$ is a probability space and $S$ is a separable Banach space with norm $\|\cdot\|$ and Borel $\sigma$-algebra $\mathscr{S}$. Let $g$ be a measurable mapping from $(\Omega,\mathscr{F})$ to $(S,\mathscr{S})$ and $M\in \mathscr{S}$ be a nonempty, closed and convex subset. Then
 \beqlb\label{3.1}
h(\omega):=\arg\min_{f\in M}\|f-g(\omega)\|
 \eeqlb
is well defined and $\mathscr{F}$-measurable.
 \end{theorem}
 \proof
 From Lemma~\ref{NPCBt2.1}, we have $h(\omega)$ exists uniquely. Now we prove it is $\mathscr{F}$-measurable.
 Let $m(\omega)=\min_{f\in M}\|f-g(\omega)\|$ which is $\mathscr{F}$-measurable. Indeed, since $M$ is separable, there exists a countable subset of $M$, denoted by $M_1:=\{f_1,f_2,\cdots\}$ such that
 $$m(\omega)=\lim_{n\rightarrow\infty}\min_{1\leq i\leq n}\|f_i-g(\omega)\|.$$
 Since $\|f_i-g(\omega)\|$ and $\min_{1\leq i\leq n}\|f_i-g(\omega)\|$ are measurable for any $i,n\geq 1$,
 we have  $m(\omega)$ is $\mathscr{F}$-measurable.
 Let $\mathscr{M}:=\{M\cap B: B\in \mathscr{S}\}$ and
  $\mathscr{G}=\{A\in \mathscr{M}: h^{-1}(A)\in \mathscr{F}\}$, both of them are $\sigma$-algebras.
 We just prove $\mathscr{G}$ is a $\sigma$-algebra. Obviously, $\emptyset, M\in \mathscr{G}$.
  If $A\in\mathscr{G}$, then
  \beqnn
 h^{-1}(A^c)\ar=\ar \bigcup_{f\in A^c}\{\omega: \|f-g(\omega)\|=m(\omega)\}\cr
 \ar=\ar
 \Big(\bigcup_{f\in A}\{\omega: \|f-g(\omega)\|=m(\omega)\}\Big)^c\in\mathscr{F}
  \eeqnn
   For any $\{A_n\}_{n=1}^{\infty}\in \mathscr{G}$, then
  \beqnn
  h^{-1}\Big(\bigcup_{n=1}^{\infty}A_n\Big)\ar=\ar\bigcup_{f\in \bigcup_{n=1}^{\infty}A_n}\{\omega: \|f-g(\omega)\|=m(\omega)\}\cr
  \ar=\ar\bigcup_{n=1}^{\infty}\bigcup_{f\in A_n}\{\omega: \|f-g(\omega)\|=m(\omega)\}\in\mathscr{F}.
  \eeqnn
 For any $a>0$ and $\xi\in M$, let
 $A=\{f\in M: \|f-\xi\|\leq a \}\in\mathscr{M}.$ The desired result follows if we prove for any  $A\in \mathscr{G}$ have
 \beqnn
 h^{-1}(A)=\bigcup_{f\in A}\{\omega: \|f-g(\omega)\|=m(\omega)\}\in \mathscr{F}.
 \eeqnn
 Since $M$ is separable, for any $i>0$ there exits a subset $\{f_1^i,f_2^i,\dots\}\subset A$ such that
 for any $f\in A$ there exists an element $f_{n}^i$ with $\|f_{n}^i-f\|<1/i$ .
 Let
 \beqnn
 B=\bigcap_{i=1}^{\infty}\bigcup_{n=1}^{\infty}\left\{\omega: m(\omega)\leq \|f_n^i-g(\omega)\|\leq m(\omega)+1/i\right\}\in\mathscr{F}.
 \eeqnn
 So it suffices to prove $h^{-1}(A)=B$. Actually, for any $\omega\in h^{-1}(A)$, there exists $f(\omega)\in A$ such that
 $\|f(\omega)-g(\omega)\|= m(\omega)$.
 Then for any $i>0$, there exists $ f_n^i$ satifying $\|f_n^i-f(\omega)\|<1/i$.
 Thus
  $$m(\omega)\leq \|f_n^i-g(\omega)\|\leq \|f_n^i-f(\omega)\|+\|f(\omega)-g(\omega)\|\leq m(\omega)+1/i,$$
  which means $\omega\in B$ and $h^{-1}(A)\subset B$.
  Otherwise, for any $\omega \in B$ and $i>0$, there exists $f_{n_i}^i(\omega)$ such that
   $$ m(\omega)\leq \|f_{n_i}^i(\omega)-g(\omega)\|\leq m(\omega)+1/i,$$
   which means $\|f_{n_i}^i(\omega)-g(\omega)\|\rightarrow m(\omega)$ as $i\rightarrow \infty$.
   So there exists $f(\omega)\in M$ such that $\|f(\omega)-g(\omega)\|=m(\omega)$.
    Since $S$ is a Hilbert space and $A$ is a closed and convex ball, from Theorem~2.1 we have $\arg\min_{f\in A}\|f-g(\omega)\|$ exists uniquely. Otherwise, since
    $$\min_{f\in A}\|f-g(\omega)\|\geq \min_{f\in M}\|f-g(\omega)\|\quad \mbox{and}\quad \lim_{i\rightarrow \infty}\|f_{n_i}^i-g(\omega)\|= m(\omega),$$
  we have
  $$\|f(\omega)-g(\omega)\|=\min_{f\in A}\|f-g(\omega)\|=m(\omega).$$
   So $f\in A$ and $B\subset h^{-1}(A)$. Here we have finished this proof.
 \qed

 Define
 $$\mathcal {L}^2(w):=\Big\{f(\lambda): \|f\|_{w}^2:= \int_0^{\infty}|f(\lambda)|^2w(\lambda)d\lambda<\infty\Big\},$$
 which is a Hilbert Space with inner product
 $\langle f,g\rangle_w=\int_0^{\infty}f(\lambda)g(\lambda)w(\lambda)d\lambda$
 for any $f,g\in \mathcal {L}^2(w)$. Let $\Psi:=\mathbf{T}K=\{\mathbf{T}f(\lambda):f\in K\}$, it is easy to see $\Psi \subset\mathcal {L}^2(w)$. Indeed, for any $f(z)\in K$ we have
 \beqnn
 \mathbf{T}f(\lambda)\ar=\ar  \int_0^1\int_0^{\infty}(1-e^{-z\frac{be^{-bs}\lambda}{b+c\lambda(1-e^{-bs})}})f(z)dzds\cr
 \ar\leq\ar\int_0^1\int_0^{\infty}(1-e^{-ze^{-bs}\lambda})f(z)dzds
 \leq  \int_0^{\infty}(1-e^{-ze^{-b}\lambda}) f(z)dz.
 \eeqnn
 Since $S_w$ is compact, there exists $a>0$ such that $S_w\subset [0,\theta]$. We have
 \beqlb\label{NPCB2.2}
 \|\mathbf{T}f\|^2_{w}\ar\leq\ar \int_0^{\theta}\Big|\int_0^{\infty}(1-e^{-ze^{-b}\lambda}) f(z)dz\Big|^2d\lambda\cr
 \ar\leq\ar \theta\Big|\int_0^{\infty}(1-e^{-ze^{-b}\theta}) f(z)dz\|^2\cr
 \ar\leq\ar  C\Big|\int_0^{\infty}(1\wedge z) f(z)dz\Big|^2= C |\mu(|f|)|^2<\infty,
 \eeqlb
 where $C>0$ is a constant independent to $f(z)$.
 \begin{lemma}\label{NPCBt2.3}
 For any two probability measures $Q_1(\cdot)$ and $Q_2(\cdot)$ on $\mathbb{R}_+$, we have $Q_1\overset{\rm d}=Q_2$ if and only if $L_{Q_1}(\lambda)=L_{Q_2}(\lambda)$ on some interval $[\lambda_1,\lambda_2]$ with $0\leq \lambda_1<\lambda_2$.
 \end{lemma}
 \proof
 Sufficiency is obvious, we just need to prove necessity. Since
 $L_{Q_1}(\lambda)$ and $L_{Q_2}(\lambda)$ are analytic on this strip $\{\lambda=\lambda_1+i\lambda_2\in\mathbb{C}: \lambda_1\geq0, |\lambda_2|<1\}$. By the assumption in this lemma and theorem in Brown and Churchill (2009, p.84), we have $L_{Q_1}(\lambda)=L_{Q_2}(\lambda)$ on this strip. The desired result follows from $L_{Q_1}(\lambda_1)=L_{Q_2}(\lambda_1)$ when we choose $\lambda_2=0$.
 \qed

 \begin{lemma}\label{NPCBt2.4}
 $\mathbf{T}: K\mapsto\Psi$ is a continuous bijection.
 \end{lemma}
 \proof
  It is easy to see $\mathbf{T}$ is a one-to-one mapping.
  Indeed, for any $k_1(z),k_2(z)\in K$ satisfying
 $
 \|\mathbf{T}k_1(\lambda)-\mathbf{T}k_2(\lambda)\|_w=0,
 $ we have $ \mathbf{T}k_1(\lambda)=\mathbf{T}k_2(\lambda)$ for any $\lambda\in S_w$.
  Moreover, by the definition of CBI processes,  there exist two probabilities $Q_1(\cdot)$ and $Q_2(\cdot)$ on $\mathbb{R}_+$ such that for i=1,2
 \beqnn
 L_{Q_i}(\lambda)=\exp\Big\{-\beta\int_0^1v_s(\lambda)ds-\mathbf{T}k_i(\lambda)\Big\}.
 \eeqnn
 By lemma~\ref{NPCBt2.3}, we have $Q_1\overset{d}=Q_2$ and $k_1(z)=k_2(z)$ almost everywhere. Thus $\mu(|k_1-k_2|)=0$. Since $\mathbf{T}$ is a linear operator,
 its continuity follows directly from its boundedness which have been proved in (\ref{NPCB2.2}).
\qed

 \begin{lemma}\label{NPCBt2.5}
 For any $R\geq 1$, $K_R$ and $\mathbf{T}K_R$ are compact, convex subsets of $\mathcal {L}^1(\mu)$ and $\mathcal {L}^2(w)$ respectively.
 \end{lemma}
  \proof Since $T$ is continuous, $\mathbf{T}K_R$ is compact and its convexity follows from the convexity of $K_R$ which is obvious.  It suffices to prove $K_R$ is compact.
  Since $k_R(z)$ is integrable,
 so there must exists $\{z_j\}_{j=1}^{\infty}$,
 such that $k_R(z_j)\rightarrow0$ as $j\rightarrow \infty$. For any fixed sequence $\{k_n:n=1,2,\cdots\}$ in $K_R$, we have $k_n(z_j)\rightarrow 0$ as $j\rightarrow \infty$.
 There exist $i_0\in\mathbb{N}$ and $z_0\in [i_0,i_0+1]$ such that  $k_n(z_0)\leq k_R(z_0)<1$ for for any $n\geq 1$. Without loss of generality, we assume $i_0=1$. For any $i=1,2,\cdots$, since ${\rm V}_{i}^{i+1}k_n(z)\leq R$, we have $0\leq k_n(z)\leq iR+1$ for any $ z\in[i,i+1]$.
 Similarly,  since ${\rm V}_{(i+1)^{-1}}^{i^{-1}}k_n(z)\leq R$, we have $0\leq k_n(z)\leq iR+1$ for any $z\in[(i+1)^{-1},i^{-1}]$.
 Thus $\{k_n(z):n=1,2,\cdots\}$ are uniformly bounded and have bounded variation on
 $[(i+1)^{-1},i+1]$, which means there exit two sequences of nonnegative,
 monotone increasing and right-continuous functions $\{k_{n,1}:n=1,2,\cdots\}$
 and $\{k_{n,2}:n=1,2,\cdots\}$ such that $k_n= k_{n,1}-k_{n,2}$. Obviously, we have
 \beqnn
     {\rm V}_{(i+1)^{-1}}^{i+1}k_{n,1}(z)\leq 2(i+1)R \quad \mbox{and}\quad
     {\rm V}_{(i+1)^{-1}}^{i+1}k_{n,2}(z)\leq 2(i+1)R.
 \eeqnn
 Without loss of generality, we assume that $k_{n,2}((i+1)^{-1})=0$, so $k_{n,2}(i+1)\leq 2(i+1)R$ and
 \beqnn
      k_{n,1}((i+1)^{-1})\leq iR+1, \quad k_{n,1}(i+1)\leq (3i+2)R+1.
 \eeqnn
  Applying the diagonalization argument to $k_{n,1}(z)$ and $k_{n,2}(z)$, there exists a subsequence $\{n_j^i:j=1,2,\cdots\}$ such that $k_{n_j^i,1}$ and $k_{n_j^i,2}$ converge to some functions $k^i_1$ and $k^i_2$ on $\mathbb{Q}\cap [(i+1)^{-1},i+1]$ respectively. For any and $l=1,2$, define
   \beqnn
   k^i_l(z)=\left\{\begin{array}{ll}
   \inf\left\{k^i_l(x):z\leq x\in \mathbb{Q}\cap [(i+1)^{-1},i+1]\right\},& \mbox{if } z\in[(i+1)^{-1},i+1];\cr
   0,& \mbox{otherwise}.
   \end{array}\right.
   \eeqnn
 Define $k^i=k^i_1-k^i_2$, then $k^i\in K_R$ and $k_{n_j^i}(z)\rightarrow k^i(z)$ at all continuity points of $k^i$ as $j\rightarrow \infty$. Since $k^i_1$ and $k^i_2$ are right-continuous and nonnegative monotone increasing functions, so is $k^i$ and the number of discontinuity points of $k^i$ is at most countable. Thus $k_{n_j^i}\rightarrow k^i$ almost everywhere as $j\rightarrow \infty$.
 Applying the diagonalization argument again, for $l=1,2$ we can find a subsequence of $\{k_{n_j^i,l}:j=1,2,\cdots\}$, denoted by $\{k_{n_j^{i+1},l}:j=1,2,\cdots\}$, such that the new subsequences converges to some functions $k^{i+1}_l$ on $\mathbb{Q}\cap [(i+2)^{-1},i+2]$, then repeat the program above again and we get  $k_{n_j^{i+1}}\overset{\rm{a.s.}}\longrightarrow k^{i+1}$ on $ [(i+2)^{-1},i+2]$. Obviously, $k^i(z)=k^{i+1}(z)$ on $[(i+1)^{-1},i+1]$. Thus $\{k^i:i=1,2,\cdots\}$ converges to some function $k'(z)$ almost everywhere on $[0,\infty)$. Applying the diagonalization argument again and the dominated convergence theorem, we have $\mu(|k_{n_i^i}-k'|)\rightarrow 0$ as $i\rightarrow \infty$.  Here we have proved that $K_R$ is compact.
 \qed

  \begin{corollary}\label{NPCBt2.6}
 For any $R\geq 1$, the inverse operator of $\mathbf{T}$, $\mathbf{T}^{-1}: \mathbf{T}K_R\mapsto K_R$, is continuous.
 \end{corollary}

 \noindent\textit{Proof of Theorem~\ref{NPCBt1.1}.} Define
 \beqlb\label{NPCB2.3}
 \hat{g}_{R,n}(\lambda)=\arg\min_{g(\lambda)\in \mathbf{T}K_R} \|g_{n}-g\|_w
 \eeqlb
 From Lemma~\ref{NPCBt2.5} and Theorem~\ref{NPCBt2.2}, we have $\hat{g}_{R,n}(z)$ is well defined, i.e. it exists uniquely and is $(\mathscr{F}_n)$-measurable. The desired results follows from $\hat{k}_{R,n}(z)=\mathbf{T}^{-1}\hat{g}_{R,n}(z)$ and Corollary~\ref{NPCBt2.6}.\qed

 \section{Consistency and asymptotic risk bound}
 \setcounter{equation}{0}
 \medskip

 We will prove Theorem~\ref{NPCBt1.2} in this section. The consistency of $\hat{k}_{R,n}(z)$ comes directly from the following result.
  \begin{corollary}\label{NPCBt3.1}
 $L_{n}(\lambda)\overset{\rm a.s.}\longrightarrow e^{-\mathbf{T}k(\lambda)}$ uniformly on any compact subset, i.e. for any compact subset $A$ of $\mathbb{R}_+$, we have
 $$\sup_{\lambda\in A}|L_{n}(\lambda)-e^{-\mathbf{T}k(\lambda)}|\overset{\rm a.s.}\longrightarrow 0.$$
 \end{corollary}
 \proof Obviously, $e^{-\mathbf{T}k(\lambda)}$ is Laplace transform of the distribution $\eta_0$ of $X(1)$, where $X(t)$ is a CBI processes defined by (\ref{NPCB1.2}) with $X(0)=0$ and $\beta=0$. For any $n\geq 1$, $L_{n}(\lambda)$ is Laplace transform of some measure denoted by $\mu_n$.
  Since $L_{n}(\lambda)\overset{\rm a.s.}\longrightarrow L_{\eta_0}(\lambda)$ for any
 $\lambda\in \mathbb{R}_+$, we have $\mu_n\rightarrow\eta$ weakly. The desired result follows from Lemma~7.6 in Sato (1999), i.e. $L_{n}(\lambda)=L_{\mu_n}(\lambda)\rightarrow L_{\eta}(\lambda)$
 uniformly on any compact subset.
 \qed

 \noindent\textit{Proof of Theorem~\ref{NPCBt1.2}.} Recall $\hat{g}_{R,n}(\lambda)$ defined in the proof of Theorem~\ref{NPCBt1.1}. From Corollary~\ref{NPCBt3.1} and Theorem~7.6.3 in Chung (2001), we have
 $$\sup_{\lambda\in S_w}|L_{n}(\lambda)-e^{-\mathbf{T}k(\lambda)}|\overset{\rm a.s.}\rightarrow 0\quad \mbox{and}\quad \sup_{\lambda\in S_w}|g_{n}(\lambda)-\mathbf{T}k(\lambda)|\overset{\rm a.s.}\rightarrow 0.$$
 By the definition of $\hat{g}_{R,n}$ in (\ref{NPCB2.3}),
 \beqlb\label{NPCB3.1}
 \|\hat{g}_{R,n}-\mathbf{T}k(\lambda)\|_w\ar\leq\ar \|\hat{g}_{R,n}-g_{n}\|_w+\|g_{n}-\mathbf{T}k(\lambda)\|_w\leq  2\|g_{n}-\mathbf{T}k(\lambda)\|_w\overset{\rm a.s.}\rightarrow 0.
 \eeqlb
 The desired result follows from this result and Corollary~\ref{NPCBt2.6}.
 \qed

 \begin{lemma}\label{NPCBt3.2}
 Suppose conditions in Theorem~\ref{NPCBt1.3} hold, then there exists a constant $C>0$ such that for $n$ large enough have
 $$n\textbf{E}\big[\|\hat{g}_{R,n}-\mathbf{T}k\|^2_w\big]< C.$$
 \end{lemma}
 \proof Let $G_n(\lambda)=\sum_{k=1}^n \xi_k(\lambda)$, where
 $$\xi_k(\lambda)=\exp\Big\{-\lambda X_k+X_{k-1}v(\lambda)+\int_0^1\psi(v_s(\lambda))ds\Big\}-1.$$
 It's easily to prove $\{G_n(\lambda)\}_{n=1}^{\infty}$ is a $(\mathscr{F}_n)$-martingale.
 Obviously, we have
  \beqlb\label{NPCB3.2}
 \textbf{E}[\xi_{k}^2|\mathscr{G}_{k-1}]\ar=\ar \textbf{E}\Big[ \Big(\exp\big\{-\lambda X_k+X_{k-1}v(\lambda)+\int_0^1\psi(v_s(\lambda))ds\big\}-1\Big)^2|\mathscr{G}_{k-1}\Big]\cr
 \ar=\ar \textbf{E}\Big[\exp\big\{-2\lambda X_k+2X_{k-1}v(\lambda)+2\int_0^1\psi(v_s(\lambda))ds\big\}|\mathscr{G}_{k-1}\Big]-1\cr
 \ar=\ar \exp\Big\{-X_{k-1}(v(2\lambda)-2v(\lambda))-\int_0^1\psi(v_s(2\lambda))ds+2\int_0^1\psi(v_s(\lambda))ds\Big\}-1.\quad
  \eeqlb
 Define $V_n=\sum_{k=1}^n \textbf{E}[\xi_{k}^2|\mathscr{G}_{k-1}]$. From (\ref{NPCB3.2}) and ergodicity of $X(t)$ we have
 \beqnn
  \frac{1}{n}V_n\ar\rightarrow\ar \textbf{E}\left[\exp\Big\{-X_0(v(2\lambda)-2v(\lambda))-\int_0^1\psi(v_s(2\lambda))ds+2\int_0^1\psi(v_s(\lambda))ds\Big\}\right]-1= W(\lambda),
 \eeqnn
 where
 \beqnn
  W(\lambda)= \exp\Big\{-\int_0^{\infty}\psi(v_s(v(2\lambda)-2v(\lambda)))ds-\int_0^1\psi(v_s(2\lambda))ds+2\int_0^1\psi(v_s(\lambda))ds\Big\}-1.
 \eeqnn
 It is easy to identify that $W(\lambda)$ is continuous and bounded on any compact set. By the martingale central limitation theorem (see Durrett, 2010), we have
 \beqnn
 h_n(\lambda):=\frac{1}{\sqrt{n}}G_n(\lambda)=\frac{1}{\sqrt{n}}\sum_{k=1}^n \xi_k(\lambda)\overset{\rm d}\rightarrow N(0,W(\lambda)).
 \eeqnn
  Furthermore, by the definition of $g_{n}$
 \beqnn
 g_{n}(\lambda)-\mathbf{T}k(\lambda)\ar=\ar\ln\Big[\frac{1}{n}\sum_{k=1}^n\exp\big\{-\lambda X_k+X_{k-1}v(\lambda)+\int_0^1\psi(v_s(\lambda))ds\big\}\Big]=\ln\left[h_n(\lambda)/\sqrt{n}+1\right].
 \eeqnn
  By the Skorohod's representation theorem, see Billingsley (1999, p.70), there exist $\{\tilde{h}_n(\lambda):n=1,2,\cdots\}$ and $\tilde{h}(\lambda)$ defined on $([0,1],\mathscr{B}[0,1],m)$, where $m$ is Lebesgue measure, such that $\tilde{h}_n(\lambda)\overset{\rm d}=h_n(\lambda)$ for any $n\geq 1$, $\tilde{h}(\lambda)\overset{\rm d}=N(0,W(\lambda))$ and $\tilde{h}_n(\lambda)\overset{\rm a.s}\rightarrow \tilde{h}(\lambda)$. Define $\tilde{H}_n(\lambda)= \ln[\tilde{h}_n(\lambda)/\sqrt{n}+1]$, then
 \beqnn
 \sqrt{n}(g_{n}(\lambda)-\mathbf{T}k(\lambda))\overset{\rm d}=\sqrt{n}\tilde{H}_n(\lambda)\overset{\rm a.s.}\longrightarrow \lim_{n\rightarrow \infty}\tilde{h}_n(\lambda)\overset{\rm d}= N(0,W(\lambda)),
 \eeqnn
 thus $\sqrt{n}(g_{n}(\lambda)-\mathbf{T}k(\lambda))\overset{\rm d}\rightarrow N(0,W(\lambda))$ as $n\rightarrow \infty$ and
 \beqnn
 \lim_{n\rightarrow \infty}\mathbf{E}[n|g_{n}(\lambda)-\mathbf{T}k(\lambda)|^2]=W(\lambda).
 \eeqnn
 By Tonelli's theorem and Fatou's theorem, we have
 \beqnn
 \limsup_{n\rightarrow\infty}n\textbf{E}[\|g_{n}-\mathbf{T}k\|^2_w]\ar=\ar \limsup_{n\rightarrow\infty}\int_0^{\infty}\textbf{E}\big[n|g_{n}(\lambda)-\mathbf{T}k(\lambda)|^2 \big]w(\lambda)d\lambda\cr
 \ar\leq\ar  \int_0^{\infty} \limsup_{n\rightarrow\infty}\textbf{E}\left[n|g_{n}(\lambda)-\mathbf{T}k(\lambda)|^2\right]w(\lambda)d\lambda\cr
 \ar=\ar\int_0^{\infty}W(\lambda)w(\lambda)d\lambda<\infty.
 \eeqnn
 From this and(\ref{NPCB3.1}), there exists a constant $C>0$ such that
 $$n\textbf{E}[\|\hat{g}_{R,n}-\mathbf{T}k\|_w^2]\leq 4n\textbf{E}[\|g_{n}-\mathbf{T}k\|^2_w]<C.$$
 Here we have finished this proof.
 \qed

 \noindent\textit{Proof of Theorem~\ref{NPCBt1.3}.}
 Obviously, we have
 \beqnn
 \mu(|\hat{k}_{R,n}-k|)= \mu(|\mathbf{T}^{-1}\hat{g}_{R,n}-k|)=\mu(|\mathbf{T}^{-1}(\hat{g}_{R,n}-\mathbf{T}k)|).
 \eeqnn
 For any fixed $R\geq 1$, since $\mathbf{T}^{-1}$ is a continuous and linear bijection from $\mathbf{T}K_R$ to $K_R$, from Corollary~2.12(c) in Rudin (1991, p49) we have
 \beqnn
 \mu(|\mathbf{T}^{-1}(\hat{g}_{R,n}-\mathbf{T}k)|)<C \|\hat{g}_{R,n}-\mathbf{T}k\|_w,
 \eeqnn
 where  $C>0$ is a constant determined by the norm of $\mathbf{T}^{-1}$.
 The desired result follows from this and Lemma~\ref{NPCBt3.2}.
 \qed

 \noindent\textbf{Acknowledgements.}  Here, I would like to thank
Professor Zenghu Li for suggesting me considering this problem.

\bigskip


\begin{thebibliography}{99}

\bibitem{AT88} Ahn, C.M., and Thompson, H.E. (1988): Jump-diffusion processes and the term structure of interest rates. \textit{J. Finance}, \textbf{43(1)}, 155--174.

 \bibitem{B99}Billingsley, P. (1999): \textit{Convergence of Probability Measures}, second edition. John Wiley \& Sons.

\bibitem{BC09} Brown, J.W. and Churchill, R.V.(2009): \textit{Complex Variables and Applications}, eighth edition. New York. McGraw-Hill.

\bibitem{C01} Chung, K.L. (2001): \textit{A Course in Probability Theory}, third editon. Academic press.

\bibitem{CIR85} Cox, J., Ingersoll, J. and Ross, S. (1985): A theory
    of the term structure of interest rate. \textit{Econometrica},
    \textbf{53}, 385--408.

\bibitem{CG09} Comte, F. and Genon-Catalot, V. (2009): Nonparametric
    estimation for pure jump L\'{e}vy processes based on high
    frequency data. \textit{ Stochastic Process. Appl.},
    \textbf{119}, 4088--4123.

\bibitem{CG10} Comte, F. and Genon-Catalot, V. (2010): Nonparametric
    adaptive estimation for pure jump L\'{e}vy processes. \textit{
    Ann. Inst. H. Poincar\'{e} Probab. Statist.},
    \textbf{46}, 595--617.

\bibitem{DPS00} Duffie, D., Pan, J. and Singleton, K. (2000): Transform analysis and asset pricing for affine jump-diffusions. \textit{Econometrica}, \textbf{68(6)}, 1343--1376.

\bibitem{DFS03} Duffie, D., Filipovi\'{c}, D. and Schachermayer, W.
    (2003): Affine processes and applications in finance. \textit{
    Ann. Appl. Probab.}, \textbf{13}, 984--1053.

\bibitem{Dur10} Durrett, R. (2010): \textit{Probability: Theory and
    Examples}. Cambridge University Press.

\bibitem{F09}  Figueroa-L\'opez, J.E. (2009): Nonparametric estimation for L\'evy models based on discrete-sampling. \textit{Lecture Notes-Monograph Series}, 117--146.

\bibitem{FH06} Figueroa-L\'{o}pez, J. and Houdr\'{e}, C. (2006): Risk
    bounds for the nonparametric estimation of L\'{e}vy processes.
    In\textit{ High Dimensional Probab. IMS Lecture Notes},
    \textbf{51}, 96--116.

\bibitem{FL10} Fu, Z., and Li, Z. (2010): Stochastic equations of
    non-negative processes with jumps. \textit{ Stochastic
    Process. Appl.}, \textbf{120(3)}, 306--330.


\bibitem{HMZ11} Huang, J., Ma, C. and Zhu, C. (2011): Estimation for
    discretely observed continuous state branching processes with
    immigration. \textit{Statist.  Probab. Letters}, \textbf{81},
    1104--1111.


\bibitem{JMV05}Jongbloed, G., van der Meulen, F.H. and van der Vaart,
    A.W. (2005): Nonparametric inference for L\'{e}vy-driven
    Ornstein-Uhlenbeck processes. \textit{Bernoulli}, \textbf{11},
    759--791.

\bibitem{S99} Sato, K. (1999): \textit{L\'evy processes and infinitely divisible distributions.} Cambridge university press.


\bibitem{Li11} Li, Z. (2011): \textit{Measure-Valued Branching Markov
    Processes}. Springer, Berlin.

\bibitem{LiMa13} Li, Z. and Ma, C. (2015): Asymptotic properties of estimators in a stable Cox–Ingersoll–Ross model. \textit{Stochastic Process. Appl.}, \textbf{125(8)}, 3196--3233.


\bibitem{NR09} Neumann, M. and Reiss, M. (2009): Nonparametric estimation for L¨¦vy processes from low-frequency observations. \textit{ Bernoulli}, \textbf{15}  223--248.


\bibitem{O98} Overbeck, L. (1998): Estimation for continuous branching processes. \textit{Scand. J. Statist.}, \textbf{25(1)}, 111--126.

\bibitem{OvR97} Overbeck, L. and Ryd\'en, T. (1997): Estimation in
    the Cox-Ingersoll-Ross model. \textit{Econometric Theory},
    \textbf{13}, 430--461.

\bibitem{R91}  Rudin, W. (1991): \textit{Functional analysis}, second edition. McGraw-Hill.

\bibitem{WL03} Watteel, R.N. and Kulperger, R.J. (2003): Nonparametric
    estimation of the canonical measure for infinitely divisible
    distributions. \textit{J. Statist. Comput. Simul.}, \textbf{73},
    525--542.

\bibitem{Xu13} Xu, W. (2013): Parameter estimation in two-type continuous-state branching processes with immigration. \textit{Statist.  Probab. Letters}, \textbf{91}, 124--134.
\end{thebibliography}
\end{document}